\RequirePackage{fix-cm}

\documentclass[twocolumn]{svjour3}          % twocolumn
\smartqed  % flush right qed marks, e.g. at end of proof
\usepackage{graphicx}
 \usepackage{mathptmx}      % use Times fonts if available on your TeX system
\usepackage{amsmath}
\usepackage{amsfonts}
\usepackage{amssymb}

\usepackage{algorithm}
\usepackage{algorithmic}

\usepackage{url}

\begin{document}

\title{Modelling Character Motions on Infinite-Dimensional Manifolds%\thanks{Grants or other notes
}


\author{Markus Eslitzbichler}


\institute{M. Eslitzbichler
    \at Department of Mathematical Sciences,
    Norwegian University of Science and Technology,
    7491 Trondheim, Norway
    \\Tel.: +47 735 91650
    \\\email{markuses@math.ntnu.no}
}

\date{Received: date / Accepted: date}

\maketitle

\begin{abstract}
%abstract.tex
%Author: Markus Eslitzbichler
%Paper: Modelling Character Motions on Infinite-Dimensional Manifolds
%Submitted to: The Visual Computer

In this article, we will formulate a mathematical framework that allows us to treat character animations as points on infinite dimensional Hilbert manifolds.
Constructing geodesic paths between animations on those manifolds allows us to derive a distance function to measure similarities of different motions.
This approach is derived from the field of geometric shape analysis, where such formalisms have been used to facilitate object recognition tasks.

Analogously to the idea of shape spaces, we construct motion spaces consisting of equivalence classes of animations under reparametrizations.
Especially cyclic motions can be represented elegantly in this framework.

We demonstrate the suitability of this approach in multiple applications in the field of computer animation.
First, we show how visual artifacts in cyclic animations can be removed by applying a computationally efficient manifold projection method.
We next highlight how geodesic paths can be used to calculate interpolations between various animations in a computationally stable way.
Finally, we show how the same mathematical framework can be used to perform cluster analysis on large motion capture databases, which can be used for or as part of motion retrieval problems.

%Include keywords, PACS and mathematical
%subject classification numbers as needed.
\end{abstract}

%introduction.tex
%Author: Markus Eslitzbichler
%Paper: Modelling Character Motions on Infinite-Dimensional Manifolds
%Submitted to: The Visual Computer

\section{Introduction}\label{Sec:Intro}
Skeletal animation is a cornerstone of modern computer graphics, used in movie special effects, tv series and video games.
Even with the major advances in graphical fidelity achieved in the last few years, creating convincing animations, particularly of human characters, remains a challenge.
\emph{Motion capturing}, where an actor's motions are recorded and superimposed on a virtual model, remains one of the most popular ways of creating realistic animations.
However, it is an inherently static process, requiring extensive algorithmic work to adapt it to interactive applications.
Furthermore, as motion capture databases grow bigger, methods of efficiently querying those databases become a concern.
Given this background, we believe there to be unexplored potential to use methods from \emph{geometric shape analysis} in the field of computer animation.

Tasks in shape analysis, such as comparing different curves to measure their similarity (given a suitable definition of ``similar''), can be elegantly formulated using tools from differential geometry.
For example, a popular approach to comparing two planar curves, consists of first gathering all curves satisfying certain criteria\footnote{Such as arclength parametrization, open/closed curves etc.} in a Riemannian manifold and then finding minimal geodesics between the two given curves on this manifold \cite{klassen_analysis_2004}.
Measuring the length of such a geodesic can then give a value indicating the similarity of the two curves.
If we think of the 2 dimensional outline of an object as that object's ``shape'', represented by a planar curve, we can answer simple classification tasks, such as ``What class of object does this shape most likely represent?''.
Similar constructions can be considered for surfaces representing 3D shapes \cite{kurtek_novel_2010}.
Such methods can be further extended and generalized to deal with tasks in, e.g., protein structure comparison \cite{liu_mathematical_2011}.

In this paper, we want to highlight how similar techniques can be advantageously applied in the field of computer animation, and more specifically skeletal animation.
By considering a character's motion as a high-dimensional curve in $\mathbb{T}^n$, the $n$-dimensional torus, we can bring to bear the machinery of differential geometry to computer animation.
This formalism can then be used for tasks as varied as motion synthesis, motion interpolation as well as motion recognition (or retrieval).

Even though we will only consider examples of human locomotion, the methods discussed here can be applied to much more diverse problems.

We will proceed as follows:
In Section \ref{Sec:MathFormulation}, we formulate the mathematical framework necessary to treat animations as points on manifolds and review some basic notions from shape analysis that we will use later on.

In Sect.~\ref{Sec:Applications}, we investigate multiple applications: \begin{itemize}
    \item In Sect.~\ref{Sec:CycleClosing}, we present a manifold projection algorithm to compute cyclic approximations of existing non-cyclic animations. (Cyclic animations are periodic, i.e., they can be repeated continuously, such as for example a walking motion.)
    \item Motion blending, i.e., weighted combinations between two animations, is often achieved by linear or spherical linear interpolations for translations or rotations respectively.
        This is often used to generate smooth transitions from one animation to another or to generate a combination of two different animations. In Sect.~\ref{Sec:Blending}, we demonstrate how geodesics on animation spaces can be used for such tasks.
    \item Using a geodesic distances, we show in Sect.~\ref{Sec:Classification} how the manifold modelling approach can be used to classify animations robustly.
\end{itemize}

\section{Previous Work}\label{Sec:PrevWork}
The mathematical techniques used in this work are largely derived from the field of shape-analysis of curves.
Most of the algorithms and mathematical formulations we use in the following sections are derived from \cite{bauer_overview_2014} and \cite{srivastava_shape_2011} and references therein.
For more information on the use of differential geometry in shape and image analysis, we refer to the overview papers \cite{srivastava_advances_2012} and \cite{younes_spaces_2012}.

The modelling of motions on manifolds has previously been investigated in \cite{abdelkader_silhouette-based_2011}, where human gestures were treated as a series of outlines of characters. The idea of a shape space for entire triangular meshes was used in \cite{kilian_geometric_2007} in tasks such as shape exploration, shape deformation and deformation transfer.

The problem of modelling cyclic animations has been covered in, among others, \cite{ormoneit_representing_2005} and \cite{gonzalez_castro_cyclic_2010}.
The authors of \cite{ormoneit_representing_2005} use a timeseries representation and Fourier analysis techniques, while the authors of \cite{gonzalez_castro_cyclic_2010} apply PDE methods to generate cyclic animations for humans as well as fish locomotion.
In Sect.~\ref{Sec:CycleClosing}, we will present a geometric way of working with such motions.

Motion blending has long been used to combine multiple animations in various ways, either to create transitions from one animation to another or to create a simple interpolation between two distinct but similar motions.
Typical challenges include time synchronization and root alignment.
A detailed treatment of this topic can be found in \cite{kovar_flexible_2003}.
Note especially the use of ``timewarp curves'' to synchronize two different animations, somewhat similar to the reparametrizations we discuss in the next section.

Tackling problems in motion retrieval has spawned a large body of work, using techniques such as time-series analysis to compare motions, indexing and search algorithms as well as work on dimensionality reduction for motions.
We refer to the review paper \cite{pejsa_state_2010} for an overview of this field of research.
%math_formulation.tex
%Author: Markus Eslitzbichler
%Paper: Modelling Character Motions on Infinite-Dimensional Manifolds
%Submitted to: The Visual Computer

\section{Mathematical Formulation}\label{Sec:MathFormulation}
\subsection{Skeletal Animation}\label{Sec:SkeletalAnimation}
\begin{figure}
    \center
    %\vspace{-1cm}
    %\includegraphics[scale=0.4]{skeleton.pdf}
    \includegraphics[scale=0.4]{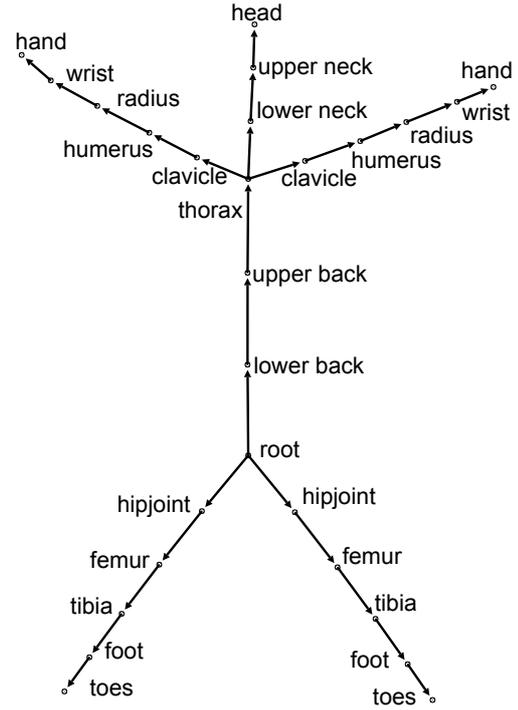}
    \caption{\label{Fig:Skeleton16} The human skeleton used for computer animation in our examples. Based on the \emph{CMU Graphics Lab Motion Capture Database} \cite{carnegie-melllon_carnegie-mellon_2003}. Note that this is not to scale}
\end{figure}

When we talk about animations, we mean \emph{skeletal animation}, where the motions of characters are defined by a skeleton and an animation curve.
The skeleton, just like the name suggests, is a hierarchy of \emph{bones} connected by \emph{joints} that defines transformation relationships.
This hierarchy is represented by a directed acyclic graph in which every node has at most one parent.
Figure \ref{Fig:Skeleton16} shows the skeleton we will be using in our example applications in Section \ref{Sec:Applications}.
Every bone is connected to its parent by a joint, which represents a transformation with respect to the parent bone.
By traversing the graph starting at a \emph{root bone} and composing all transformations along the path, we get a global transformation for every bone. \footnote{The root bone can also have rotational and translational transformations with respect to a global coordinate system.}

The transformation from a parent to a child bone is analogous to a joint in a real skeleton.
Note however, that the skeletons used in computer animation are generally much simpler than real life skeletons.
They are designed to replicate or model a given range of human motions, while keeping in mind computational performance.
Therefore bones or joints are often left out, grouped together, or new, artificial bones are introduced.

Every bone has a fixed length\footnote{i.e., translation w.r.t to its parent.} and there are one to three degrees of rotational freedom associated with it.\footnote{By a bone's degrees of freedom we mean, more precisely, the corresponding joint's degrees of freedom.}
Taking for example a person's right leg, we see that the tibia (shinbone) is a child bone of the femur (thigh bone), with a single degree of rotation between them - the knee.
The tibia, in turn, is the parent bone to the ``foot bone'' which has two degrees of freedom - pitch and yaw.

Furthermore, degrees of freedom can have constraints placed on them - we don't want a knee to bend backwards for example.
But for now, we will ignore such constraints.

We can collect all bones in the set $\mathcal{B}$ and denote a bone's number of degrees of freedom by $\text{dof}(b)$ for $b \in \mathcal{B}$.
We can then define \emph{joint space} $\mathcal{J}$ as: \begin{displaymath}
    \mathcal{J} := \mathbb{T}^n = \underbrace{\mathbb{S}^1 \times \cdots \times \mathbb{S}^1}_{n},
\end{displaymath} where $n = \sum_{i \in \mathcal{B}} \text{dof}(i)$ is the total number of degrees of freedom in the skeleton, $\mathbb{T}^n$ is the $n$-dimensional torus and $\mathbb{S}^1$ denotes the unit circle.
Any given pose of the skeleton then corresponds to a single point in joint space $\mathcal{J}$.
We have chosen this Euler angle representation in order to simplify computations in the numerical examples in Sect.~\ref{Sec:Applications}.
Future work could include using elements in the special orthogonal group $\text{SO}(3)$ to represent every bone.

An \emph{animation curve} $\beta$ is a function that assigns, for every point in a given time interval, a transformation to every bone in a skeleton: \begin{displaymath}
    \beta : [0,T] \rightarrow \mathcal{J}.
\end{displaymath}

In practice, an animation is often only specified by joint angles at discrete intervals of time, in between which interpolation is used to determine character poses.
This is for example the case in motion captured animations, where an actor's motions are recorded at a given sampling rate.

While joint space $\mathcal{J}$ is well suited to describing animations, rendering character models and physics calculations such as collision detection require positions in world space, i.e., $\mathbb{R}^3$.
We denote the map computing a given bone's position in the root's coordinate system by \begin{equation}\label{Eq:JointToWorld}
    \mathcal{W} : \mathcal{B} \times \mathcal{J} \rightarrow \mathbb{R}^3.
\end{equation}
In our particular case, this is just the multiplication of successive transformation matrices along the bone hierarchy.
The final 3d world-space position is then determined by adding translation and rotation of the root bone.

%\todo{scale down figure}

\subsection{Manifolds of curves}\label{Sec:ManifoldsOfCurves}
In this section, we will review methods from curve analysis, based on \cite{bauer_overview_2014}, \cite{bauer_sobolev_2011} and \cite{srivastava_shape_2011}.
We assume that our curves are immersions, i.e., smooth maps which have injective derivatives everywhere.
For curves, this means having a non-vanishing derivative, so we can write the space as \begin{displaymath}
    \text{Imm}(M, \mathbb{R}^n) := \lbrace q \in C^\infty(M, \mathbb{R}^n): \dot{q}(t) \neq 0 \quad \forall t \in M \rbrace ,
\end{displaymath} where $M$ is either the interval $[0,T]$ for open curves or the unit circle $\mathbb{S}^1$ for closed curves.

Given a unit length curve $c \in \text{Imm}(M, \mathbb{R}^n)$, we define the mapping \begin{equation}\label{Eq:RMap}
    R[c] := \dfrac{\dot{c}}{\sqrt{\|\dot{c}\|}},
\end{equation} and we refer to $R[c]$ as the representation of $c$.
The rationale behind this mapping will become clear soon.
The inverse of $R$ is defined up to a translation: \begin{equation}\label{Eq:RMapInv}
    R^{-1}[\beta](t) := c_0 + \int_0^t \beta(s) \| \beta(s) \| ds.
\end{equation}
The following constructions can be extended to curves of arbitrary length with only minor modifications.

The image of unit length immersions under $R$ coincides with \begin{equation}\label{Eq:COpen}
    \mathcal{C}^o := \{ q \in C^\infty([0,T], \mathbb{R}^n) : \int_0^T \| q(t) \|^2 dt = 1\}.
\end{equation}
The image of closed unit length immersion under $R$ coincides with \cite{srivastava_shape_2011} \begin{align}\label{Eq:CClosed}
    \mathcal{C}^c := \{ q \in C^\infty(\mathbb{S}^1, \mathbb{R}^n) : & \int_{\mathbb{S}^1} \| q(t) \|^2 dt = 1\\
    & \wedge \int_{\mathbb{S}^1} q(t) \| q(t) \| dt = 0\}. \nonumber
\end{align}

%Let us now take a look at the geometry of those two manifolds.
Tangent vectors to points in the manifolds $\mathcal{C}^o$ and $\mathcal{C}^c$ are vector fields along the corresponding curves.
They can be interpreted as infinitesimal deformations of the curves.
For the manifold of unit length open curves, elements in the tangent space need to preserve the unit length:

\begin{equation*}
    \mathrm{T}_q\mathcal{C}^o := \{ v \in C^\infty([0,T], \mathbb{R}^n) : \ \langle v, q \rangle_{L^2([0,T], \mathbb{R}^n)} = 0\}.
\end{equation*}

For the manifold of closed curves, tangent vectors become a bit more complicated.
The easiest way to describe the tangent space to $\mathcal{C}^c$ is via its normal space \cite{srivastava_shape_2011}:
\begin{align}
    \mathrm{N}_q\mathcal{C}^c := &\  \text{span} \lbrace q, (\dfrac{q_i}{\|q\|} q + \|q\| e_i)_{i=1,\ldots,n}\rbrace, \label{Eq:NormalSpace}
\end{align} where $e_i$ denotes the $i$'th unit basis vector in $\mathbb{R}^n$ and $q_i$ is the $i$'th component of the vector valued function $q$, and so
\begin{align*}
    \mathrm{T}_q\mathcal{C}^c := &\  \lbrace v \in C^\infty(\mathbb{S}^1, \mathbb{R}^n) : \nonumber \\
     & \ \langle v, w \rangle = 0  \qquad \forall w \in \mathrm{N}_q\mathcal{C}^c \rbrace. \nonumber
\end{align*}

We will use the standard $L^2([0,T], \mathbb{R}^n)$ inner product, restricted to $\mathcal{C}^o$ and $\mathcal{C}^c$ respectively, as a metric.
Note that this metric will, when pulled back to $\text{Imm}(M, \mathbb{R}^n)$ via the mapping (\ref{Eq:RMap}), result in a Sobolev type metric\footnote{Also known as elastic metric.} \cite{bauer_overview_2014}.
Such metrics can be used to measure bending and stretching deformations in a curve.
Intuitively, this point of view arises by thinking of a curve as an elastic string and considering the changes in elastic energy in the string under such mechanical deformations.
These Sobolev type metrics\footnote{They are parametrized by weights balancing the influences of bending vs. stretching forces.} are popular choices for shape analysis applications, but depending on parameter choices and order of the metric, they can be computationally expensive to use.
The advantage then, of representing curves using the mapping (\ref{Eq:RMap}), is that we can perform calculations using the simpler, location-invariant $L^2$ inner product, instead of the more complicated elastic metric.

Given this geometric setup, we can calculate geodesic paths on $\mathcal{C}^o$ and $\mathcal{C}^c$.
Since $\mathcal{C}^o$ is geometrically a sphere (see Definition (\ref{Eq:COpen})), there exists a simple explicit expression for the geodesic between $\beta, \gamma \in \mathcal{C}^o$ \cite{srivastava_shape_2011}: \begin{align}
    &\alpha : [0,1] \rightarrow \mathcal{C}^o \nonumber\\
    &\alpha(\tau)(t) := \dfrac{1}{\sin(\theta)} ( \sin(\theta(1-\tau)) \beta(t) + \sin(\theta \tau) \gamma(t)) \nonumber,
\end{align} where $\theta = \cos^{-1}( \langle \beta, \gamma \rangle)$.
This is simply spherical linear interpolation \cite{shoemake_animating_1985}.

For the space of representations of closed curves, $\mathcal{C}^c$, we need a numerical method to calculate the geodesic distance between $\beta, \gamma \in \mathcal{C}^c$.
Two different approaches can typically be used.
Whereas the authors of \cite{srivastava_shape_2011} use a gradient descent method which they call \emph{path-straightening}, an ODE based shooting method can also be employed to solve the geodesic boundary value problem - see \cite{bauer_overview_2014} and references therein.

For our applications we employ the gradient descent method.
We refer to \cite{srivastava_shape_2011}, Section 4.3, for details on the algorithm.

%Note that the advantage of this framework is that it allows us to use an elastic metric for curves but still perform all computations using the normal $L^2$ inner product, which is location-invariant and therefore efficient to compute.

Finally, given a geodesic path $\alpha$ from $\beta$ to $\gamma$ on either $\mathcal{C}^o$ or $\mathcal{C}^c$, we can compute the length of the geodesic via the functional \begin{equation}\label{Eq:LengthFunctional}
    L[\alpha] := \int_0^1 \| \dfrac{d}{d\tau}\alpha(\tau) \| d\tau,
\end{equation} by noting that, in a sufficiently small neighborhood, $L[\alpha]$ is minimized by a unique geodesic.

We will refer to the value of this integral as the geodesic distance between $\beta$ and $\gamma$.

\subsubsection{Shape spaces of closed curves}\label{Sec:ShapeSpaces}

\begin{figure*}
    \center
    \includegraphics[width=\textwidth]{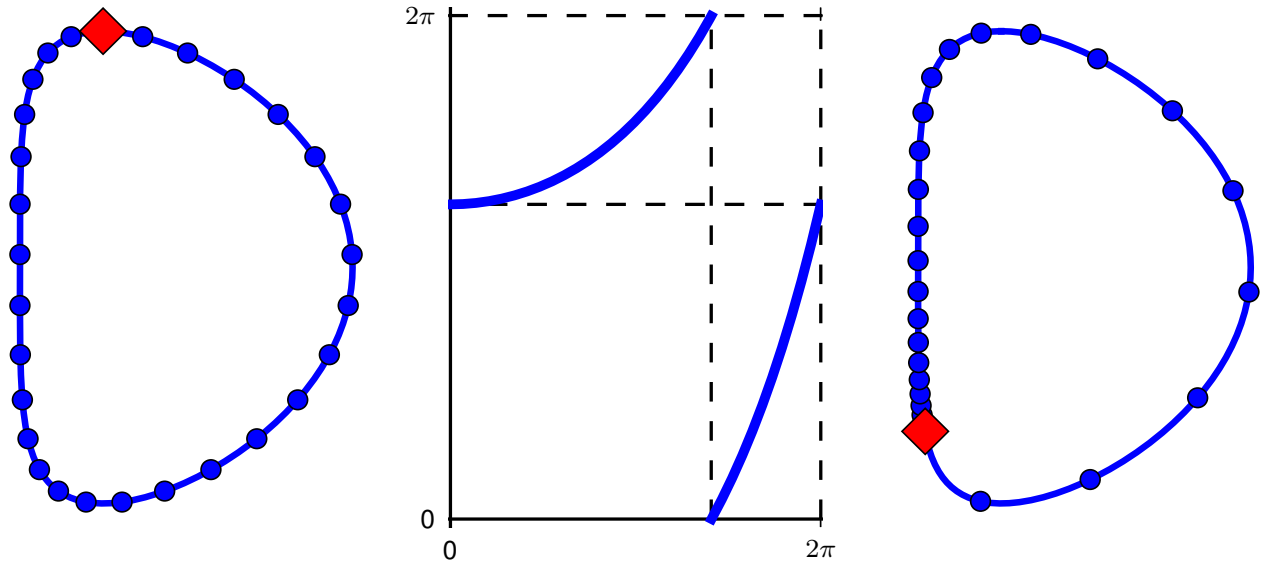}
    \caption{\label{Fig:DiffeoExample} Reparametrization of a curve by composition with a diffeomorphism on the unit circle $\mathbb{S}^1$. Both curves evolve clockwise. Note the changed starting point (the diamond) and increased density of sampling points close to the start}
\end{figure*}

If we look at the manifold of closed curves (\ref{Eq:CClosed}), we see that curves that have the same image, but different parametrizations\footnote{See Fig. \ref{Fig:DiffeoExample}}, are considered distinct in $\mathcal{C}^c$, i.e., there is no \emph{reparametrization}-invariance.
In many applications, this is undesirable.
For example in shape analysis, a circle \begin{displaymath}
    t \rightarrow (\cos(t + p), \sin(t + p)),
\end{displaymath}  is a circle, no matter what value we pick for the phase $p$.
Similarly, for example when searching in a motion capture database for a cyclic walking animation, we might choose not to care if the animation starts with the left foot or the right foot moving forward.

Such symmetries can be taken into account by introducing so-called \emph{shape spaces}.
Note that we have already accounted for translational symmetries (i.e., a circle is a circle, no matter where it is located) by using the mapping (\ref{Eq:RMap}), where only first order derivatives of the original function appear.

We can model reparametrizations as diffeomorphisms on the unit circle $\mathbb{S}^1$, i.e., smooth invertible maps from $\mathbb{S}^1$ onto $\mathbb{S}^1$, which we denote by $\text{Diff}(\mathbb{S}^1)$.
The group $\text{Diff}(\mathbb{S}^1)$ acts on the manifold $\mathcal{C}^c$ from the right: \begin{displaymath}
    \mathcal{C}^c \times \text{Diff}(\mathbb{S}^1) \rightarrow \mathcal{C}^c : (\beta, \varphi) \mapsto \beta \circ \varphi .
\end{displaymath}
See Figure \ref{Fig:DiffeoExample} for an example.

We can then define the shape space $\mathcal{S}$ as the quotient space \begin{equation}\label{Eq:ShapeSpace}
    \mathcal{S} := \mathcal{C}^c/\text{Diff}(\mathbb{S}^1).
\end{equation} 
This means that points in $\mathcal{S}$ are equivalence classes of closed curves under reparametrizations, i.e., where we have factored out starting point and speed along the curve.\footnote{We could also construct a similar space over the space of open curves $\mathcal{C}^o$, where the start point remains fixed but the sampling rate along the curve varies.}
Going back to our circle analogy, this means that we consider two curves equivalent, if we can find a starting point and a speed function to match one to another.

The space $\mathcal{S}$ inherits a Riemannian structure from $\mathcal{C}^c$ which allows for a simple characterization of geodesics in $\mathcal{S}$.
We again refer to \cite{srivastava_shape_2011} and \cite{bauer_sobolev_2011} for details.
The geodesic distance between two elements $[\beta]$ and $[\gamma]$ in $\mathcal{S}$ can then be found by \begin{equation}\label{Eq:ShapeSpaceDistance}
    d_\mathcal{S}([\beta], [\gamma]) = \inf_{\varphi \in \text{Diff}(\mathbb{S}^1)} d_{\mathcal{C}^c}( \beta, \gamma \circ \varphi),
\end{equation} where $d_{\mathcal{C}^c}$ denotes the distance function on $\mathcal{C}^c$ defined by minimizing the length functional (\ref{Eq:LengthFunctional}).
We see that, not surprisingly, finding a geodesic path between the equivalence classes $[\beta]$ and $[\gamma]$ in $\mathcal{S}$ is more involved than finding a geodesic path on $\mathcal{C}^c$.

We refer to \cite{srivastava_shape_2011} and \cite{bauer_overview_2014} for two example approaches to solving this optimization problem.

In Section \ref{Sec:Classification}, we will use this framework to calculate the similarities of different animations.
As we will see, using the quotient structure as outlined in this section enables us to get useful results for a wide range of different animations.

\subsection{Manifolds of animation curves}\label{Sec:AnimManifolds}
We now apply the concepts of the previous section to high dimensional animation curves.

The general curves of Sect.~\ref{Sec:ManifoldsOfCurves} will now be replaced by animation curves.
Cyclic animations then correspond to closed curves, represented by points on $\mathcal{C}^c$.
If we look at the mapping (\ref{Eq:RMap}) in this context, we can interpret it as a scaled angular velocity function.

As we have seen, given two curves $c$ and $d$, we look for a geodesic distance under an elastic metric, which penalizes stretching and bending.
In the context of animations, this can be seen as finding a continuous deformation of one animation into another, such that the sum of changes in angular velocity and angular acceleration are minimized.
One can think of this as the most energy efficient way of changing from one motion to another.

While this is a fairly simple model, it turns out to nevertheless give convincing results in multiple applications - from the generation of animations to the classification and recognition of motions.

The more general shape space $\mathcal{S}$ (\ref{Eq:ShapeSpace}) can also be used in the computer animation setting.
In that context, we will refer to it as \emph{motion space}.
Applying reparametrizations on animations can have a big effect.
While the starting point just shifts the motion, reparametrizations can also change the speed at which parts of the animation progress, just like the sampling rate is changed along the curve in Fig. \ref{Fig:DiffeoExample}.
Such changes can for example strongly slow down some parts of an animation while speeding up others.
This way, a walking animation can be turned into a limping walk and vice versa.
Whether or not such equivalencies are desirable depends on the specific application.
For example both in blending between two different animations and in animation classification, factoring out the starting points and playback speeds of cyclic animations can be useful; see Sections \ref{Sec:Blending} and \ref{Sec:Classification}.
%applications.tex
%Author: Markus Eslitzbichler
%Paper: Modelling Character Motions on Infinite-Dimensional Manifolds
%Submitted to: The Visual Computer

\section{Applications}\label{Sec:Applications}
We will now look at a few concrete examples of how the manifold framework from Section \ref{Sec:MathFormulation} can be used in computer animation problems.
In \ref{Sec:CycleClosing}, we propose a simple strategy to improve the periodicity of an animation (e.g., reducing discontinuities or ``jerks'' when repeating a given motion).
In \ref{Sec:Blending}, we use the notion of geodesic paths as outlined above to blend between two different animations.
Finally, in \ref{Sec:Classification}, we use geodesic distances to measure similarities between animations.
This allows us to perform cluster analysis on a motion capture database.

\subsection{Cyclic Animations}\label{Sec:CycleClosing}
Assume that we want to animate a character running forward for an unspecified amount of time.
Given only a finite set of pre-recorded animations, we will need to repeat those in a periodical manner.
However, animations generated by motion capturing are limited in time and typically non-cyclic.

In order to have a continuous forward running motion, we therefore need to generate a cyclic animation based on the motion captured running animation.
Cyclic means that the animation can be played repeatedly without visible artifacts (jumps/jerks etc.) when it repeats.
While there are many ways of calculating suitable cutoff/transition points (see, e.g., \cite{kovar_motion_2002}), finding reasonably gapless motion segments in motion recordings is improbable.

Therefore, animators often need to manually adjust animations in an effort to remove artifacts when continuously looping through the same animation.
As an alternative to this, we propose to use manifold projection techniques to automate this process.

Mathematically, we are looking at the following problem:

We are given a fixed time $T$ and an animation $c : [0,T] \rightarrow \mathcal{J}$, that is almost periodic, i.e., the errors, \begin{eqnarray}
    & \| c(T) - c(0) \| \label{Eq:MinBeta}\\
    & \| \dot{c}(T) - \dot{c}(0) \| \label{Eq:MinBetaDot},%\\
    %\| t_s - t_e\| & \rightarrow \max. \label{Eq:MaxTime}
\end{eqnarray} are small but not 0.
Conditions (\ref{Eq:MinBeta}) and (\ref{Eq:MinBetaDot}) express our desire to have sufficiently smooth periodicity (in a visual sense).
As an alternative to (\ref{Eq:MinBeta}), which measures the difference in joint angles between the start and the end poses, one could consider the $\mathbb{R}^3$ world space positions of the joints to measure the similarities in character poses using the Function (\ref{Eq:JointToWorld}): \begin{equation}
   \| \mathcal{W}(b, c(T)) - \mathcal{W}(b, c(0)) \| \qquad \forall b \in \mathcal{B},\label{Eq:MinMappedBeta} 
\end{equation} where $\mathcal{B}$ denotes the set of all bones in the animation skeleton.

Our goal now is to improve the periodicity of the animation by finding an approximation $\bar{c}: [0,T] \rightarrow \mathcal{J}$  of $c$, such that: \begin{align}
    &\| \bar{c}(T) - \bar{c}(0) \| \rightarrow 0 & \label{Eq:MinBetaBar}\\
    &\| \dot{\bar{c}}(T) - \dot{\bar{c}}(0) \| \rightarrow 0. & \label{Eq:MinBetaBarDot}
\end{align}
As before, conditions (\ref{Eq:MinBetaBar}) and (\ref{Eq:MinBetaBarDot}) express our desire to avoid rapid changes in angles or angular velocities as the animation repeats.
%, while Conditions (\ref{Eq:MinBetaBarMinusHat}) and (\ref{Eq:MinBetaBarMinusHatDot}) mean that the result should still resemble the original animation.

We propose to solve this problem by formulating it in terms of the manifolds of open and closed curves, $\mathcal{C}^o$ and $\mathcal{C}^c$, respectively.
A periodic animation then corresponds to an element in $\mathcal{C}^c$, whereas a non-periodic animation will correspond to an element in $\mathcal{C}^o$.

We first map the curve $c$, specified above, into its representation under the map (\ref{Eq:RMap}): \begin{displaymath}
    \beta := R[c].
\end{displaymath}
This $\beta$ will in general be a curve in $\mathcal{C}^o$.

Remember now that $\mathcal{C}^c$ was defined by the integral constraint (\ref{Eq:CClosed}), such that for $q \in \mathcal{C}^o$:\begin{displaymath}
    q \in \mathcal{C}^c \iff \int_{\mathbb{S}^1} q(t) \| q(t) \| dt = 0.
\end{displaymath}
This means we can project $\beta$ onto $\mathcal{C}^c$ by enforcing this constraint.
In practice, this is done by iteratively minimizing the integral \cite{srivastava_shape_2011}.
We denote this operation by \begin{equation}\label{Eq:ProjectionOp}
    \mathcal{P}^c : \mathcal{C}^o \rightarrow \mathcal{C}^c.
\end{equation}

In every iteration, the projection calculates a correction vector along the normal space (\ref{Eq:NormalSpace}).
See Algorithm \ref{Alg:ProjectC} and \cite{srivastava_shape_2011} for more details.

Finally, we construct the actual animation $\bar{c}$ using the Inverse Mapping (\ref{Eq:RMapInv}): \begin{displaymath}
    \bar{c} := R^{-1}[\bar{\beta}].
\end{displaymath}

%Algorithm \ref{Alg:ProjectC} summarizes the iterative method projecting $\beta$ into $\mathcal{C}^c$.

\begin{algorithm}[h!]
\caption{\label{Alg:ProjectC} Projection from $L^2$ into $\mathcal{C}^c$. (Taken from \cite{srivastava_shape_2011}).}
\begin{algorithmic}
    \REQUIRE $q \in L^2([0,T], \mathbb{R}^n)$ \COMMENT{Curve to project}
    \REQUIRE $\{(b_i)_{i=1,\ldots,n}\}$ \COMMENT{Basis of the normal space $N_q\mathcal{C}^c$, see Definition (\ref{Eq:NormalSpace})}

    \WHILE{$\| R^{-1}[q](T) \| \geq \epsilon$}
        \STATE $J_{i,j} \gets \ \delta^i_j + 3 \int_{\mathbb{S}^1} q_i(s) q_j(s) ds, \quad i,j=1,\ldots,n$  \COMMENT{Jacobian}
        \STATE $r \gets \  R^{-1}[q](T)$ \COMMENT{Residual}
        \STATE $\beta \gets \ -J^{-1} r$ \COMMENT{Correction}
        \STATE $q \gets \ q + \sum^n_{i=1} \beta_i b_i$  \COMMENT{Update}
    \ENDWHILE
\end{algorithmic}
\end{algorithm}

\begin{figure*}
    \center
    %\vspace{-1cm}
    %\includegraphics[scale=1., trim = 0cm 1cm 0cm 0cm]{run_closing.pdf}
    \includegraphics[width=\textwidth, trim = 0cm 1cm 0cm 0cm]{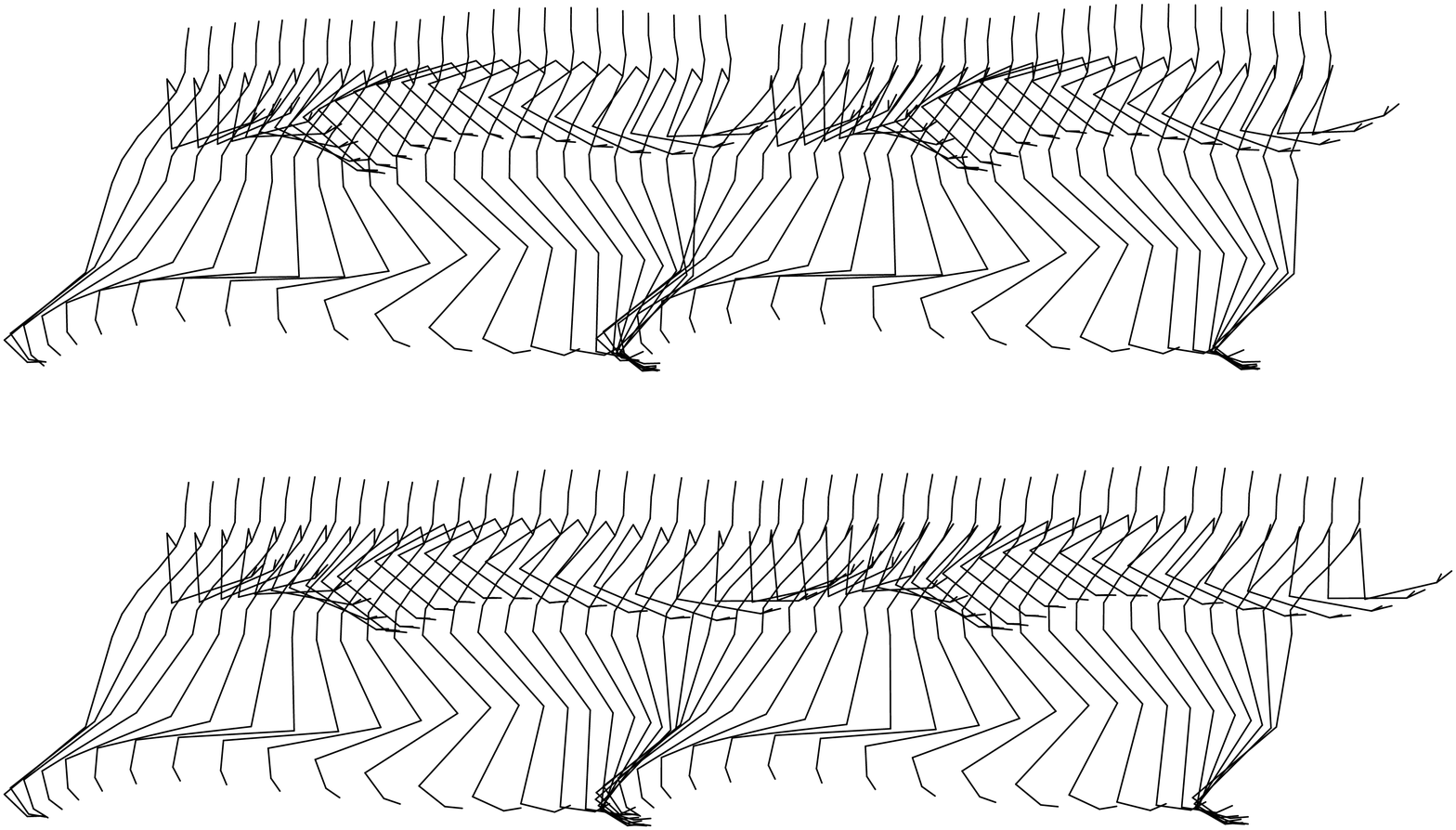}
    \caption{\label{Fig:RunCyclic} Optimization of the periodicity of a running animation. For simplicity, only half the character is shown. The top figure shows the original animation, played two times, with a noticeable gap. The bottom figure shows a more cyclic animation derived using the projection defined in Sect.~\ref{Sec:CycleClosing} and Algorithm \ref{Alg:ProjectC}}
\end{figure*}

\begin{figure*}
    \center
    %\vspace{-1cm}
    %\includegraphics[scale=0.6, trim = 2cm 1cm 0cm 4cm]{jumping_closing.pdf}
    %\includegraphics[scale=1]{jumping_closing.pdf}
    \includegraphics[width=\textwidth]{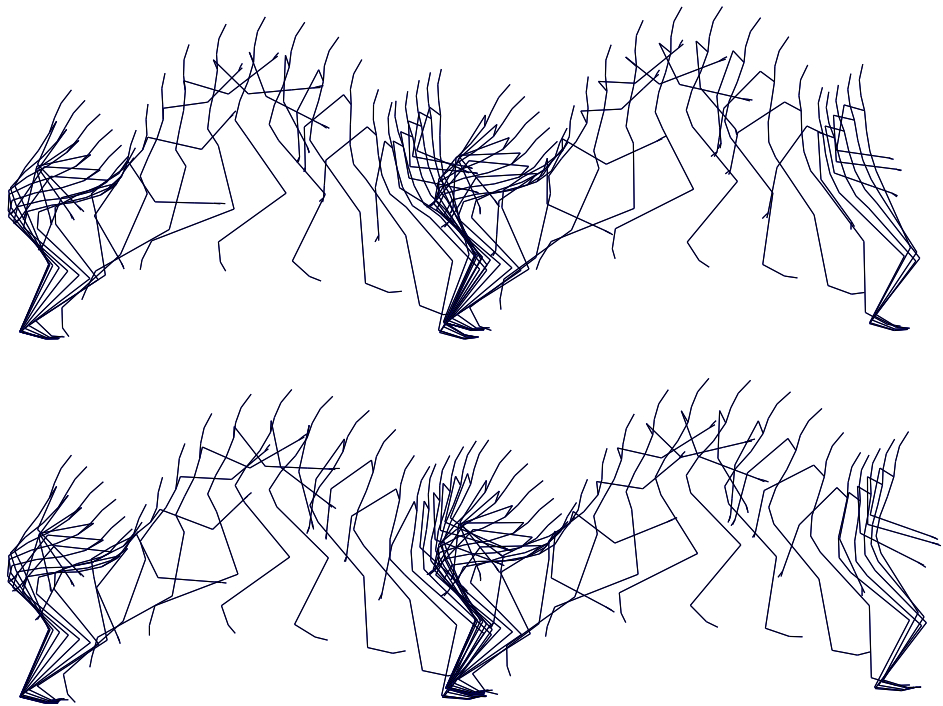}
    \caption{\label{Fig:JumpCyclic} Optimization of the periodicity of a jumping animation. For simplicity, only half the character is shown. The top figure shows the original animation, played two times, with a noticeable gap between repetitions. The bottom figure shows a more cyclic animation derived using the projection defined in Sect.~\ref{Sec:CycleClosing} and Algorithm \ref{Alg:ProjectC}}
\end{figure*}

Figures \ref{Fig:RunCyclic} and \ref{Fig:JumpCyclic} show the application of this method to running and jumping animations with noticeable jerks when repeating the motion.
We can see how the projection onto the manifold of closed curves manages to produce a much smoother repeated playback, while maintaining the original running motion without noticeable distortions.

As a further development, the use of higher-order continuity conditions to smooth transitions even more could be considered.

\subsection{Animation Blending}\label{Sec:Blending}
Motion capturing can only provide a finite set of animation samples.
In order to achieve a larger variety of animations or to continuously transition from one animation to another, various forms of \emph{motion blending} (interpolations between animations) can be used.
In its simplest form, this can be just a spherical interpolation between, for example, two walking animations to avoid too repetitive looking motions.
In order to combine more diverse animations such as walking and running, either manual modeling or more sophisticated methods are required to deal with, among others, synchronization problems \cite{kovar_flexible_2003}.

We can use geodesic paths with respect to a Riemannian metric on $\mathcal{C}^o$ and $\mathcal{C}^c$ (we will use the symbol $\mathcal{C}$ when we mean either of those two spaces) to interpolate between different animations.
Recall the notion of the geodesic path $\alpha : [0,1] \rightarrow \mathcal{C}$ between two animations $c$ and $d$ as deformations from one to the other using a minimal amount of energy.
We found that in many cases, such a path gives rise to visually convincing blends $\alpha(s)$ between $c$ and $d$ where the blend parameter $s \in [0,1]$ specifies the closeness of the interpolation to either $c$ or $d$.
This can be a simple but useful method to blend between two animations that are just slightly too different to be suitable for simple linear blending.
Note that if $\alpha(s)$ is a path on $\mathcal{C}^c$, every point along $\alpha$ corresponds to a closed curve, i.e., a cyclic animation.

As outlined in Sect.~\ref{Sec:MathFormulation}, there are various ways to calculating this geodesic path.
For our examples, we have used the ``path-straightening'' method \cite{klassen_analysis_2004}, where the set of all valid paths from $c$ to $d$ on $\mathcal{C}$ is treated as a manifold $\mathcal{H}$ itself, on which we then try to minimize the energy functional: \begin{displaymath}
    E[\alpha] = \dfrac{1}{2} \int_0^1 \langle \dfrac{d}{d\tau}\alpha(\tau), \dfrac{d}{d\tau}\alpha(\tau) \rangle d\tau,
\end{displaymath} where $\alpha: [0,1] \rightarrow \mathcal{C}$ is a path on $\mathcal{C}$ such that $\alpha(0) = c$ and $\alpha(1) = d$.
Critical points of $E$ are geodesics on $\mathcal{C}$ \cite{milnor_morse_1963}, \cite{srivastava_shape_2011}.
%By using the Palais metric \cite{palais_morse_1963}: \begin{align*}
%\langle w_0, w_1 \rangle_{T_\alpha \mathcal{H}} := & \langle w_0(0), w_1(0) >_{T_{\alpha(0)} \mathcal{C}} +\\ & \int_0^1 \langle %\frac{D w_0}{dt}(t), \frac{D w_1}{dt}(t) \rangle_{\mathbb{R}^3} dt,
%\end{align*} it becomes easy to calculate the gradient of $E$, which can then be used to efficiently determine geodesic paths on $\mathcal{C}$.
By using a special metric\footnote{The \emph{Palais}-metric; see, e.g., \cite{palais_morse_1963}.} on $\mathcal{H}$, it becomes easy to calculate the gradient of $E$, which can then be used to efficiently determine geodesic paths on $\mathcal{C}$.
We refer to \cite{klassen_analysis_2004,srivastava_shape_2011} for more details on this algorithm.

\begin{figure*}
    \center
    %\vspace{-1cm}
    %\includegraphics[scale=1., trim = 2cm 1cm 0cm 1cm]{interpolation_easy.pdf}
    %\includegraphics[scale=1.]{interpolation_easy.pdf}
    \includegraphics[width=\textwidth, height=\textheight, trim = 3cm 0cm 5cm 0cm, clip]{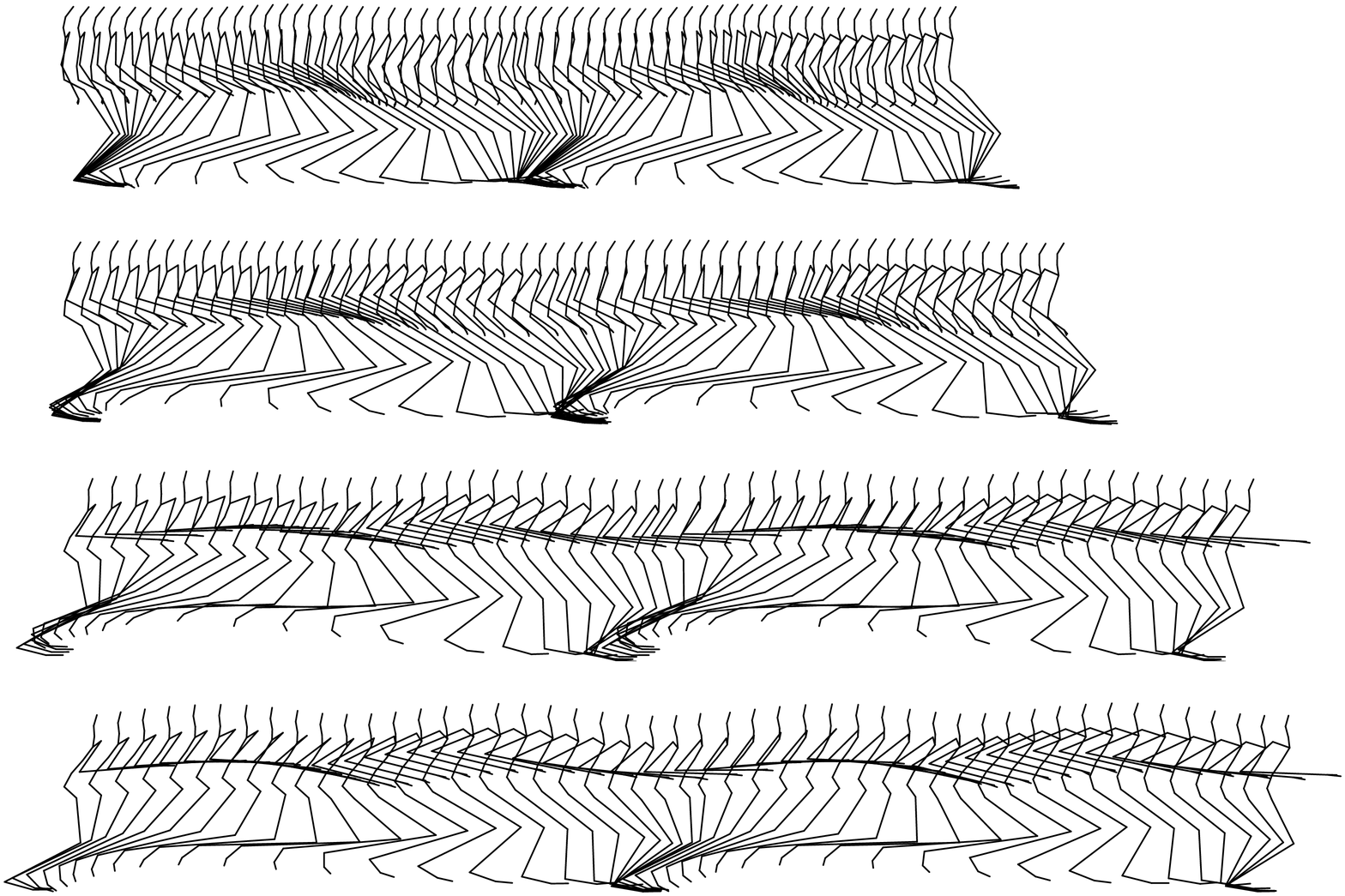}
    \caption{Interpolation between a walking (top plot) and a running animation (bottom plot) using a geodesic path. For simplicity, only half the character is shown    \label{Fig:InterpolationEasy}}
\end{figure*}

Figure \ref{Fig:InterpolationEasy} shows two such interpolated animations between a walking and a running motion.
%We can see that the energy-minimizing approach can produce visually convincing results.

For cyclic animations, this scheme can be extended by employing the animation space concept outlined in Sect.~\ref{Sec:ShapeSpaces} to calculate the geodesic path on $\mathcal{S}$ instead of $\mathcal{C}^c$. 
By factoring out the starting point and the sampling rate along a cyclic animation, a wider variety of animations can be blended together, albeit at higher computational cost.
This can be seen as somewhat analogous to \emph{timewarp curves} that are used to synchronize motions \cite{bruderlin_motion_1995,kovar_motion_2002}.

Indeed, animation curves can be annotated with additional information \cite{wei_liu_riemannian_????}.
For example, points in time where parts of a body should be stationary, such as a foot touching a ground, could be marked to improve matching performance.
This should allow for better alignment of out of sync animations, again similar to how the timewarp curves are used to synchronize multiple animations.
This will be covered in future work in this area.

Another extension to this scheme would be to use the manifold structure to create combinations of multiple animations.
As an example use, by using randomly weighted interpolations between several distinct walking animations, one could generate a continuous stream of non-repetitive forward motion.
Such weighted averages could be computed analogously to the computation of average shapes via Karcher Mean in \cite{kurtek_statistical_2012}.

\subsection{Classification}\label{Sec:Classification}
Similar to curve analysis, we can try to classify animations using the geodesic distance as defined in Sect.~\ref{Sec:MathFormulation}.
This could be used for or as part of motion retrieval applications.

Given two animations $c, d$, we map them into $\mathcal{C}$ (either $\mathcal{C}^o$ or $\mathcal{C}^c$) using (\ref{Eq:RMap}): $\beta := R[c], \ \gamma := R[d]$.
(We account for differences in the lengths of animations by rescaling along the time axis.)
We can then calculate the geodesic distance between their corresponding equivalence classes $[\beta]$ and $[\gamma]$ in the motion space $\mathcal{S}$ using Eq. (\ref{Eq:ShapeSpaceDistance}).

By calculating these distances for a large set of animations, we can construct a distance matrix, which lends itself to further statistical analysis, such as cluster identification.

As a test, we have taken the entire set of animations for ``subject 16'' in the CMU motion capture database \cite{carnegie-melllon_carnegie-mellon_2003}.
This is a set of 58 walking, jogging, running and jumping animations, with variations such as ``walk, veer left'' or ``run, 90-degree right turn'' performed by the same actor.
Every pair of animations was compared by first scaling them to a common (time) length and then calculating the geodesic distance between the two animation curves.
As a comparison test, the linear normed distance between animations was also calculated.

\begin{figure*}
    \vspace{-0.25cm}
    \begin{minipage}[l]{.6\textwidth}
        %\centering
        %\includegraphics[scale=0.8, angle=-90]{cluster16_linear.pdf}
        \includegraphics[scale=0.8, angle=-90, trim={0cm 0 1cm 0}, clip]{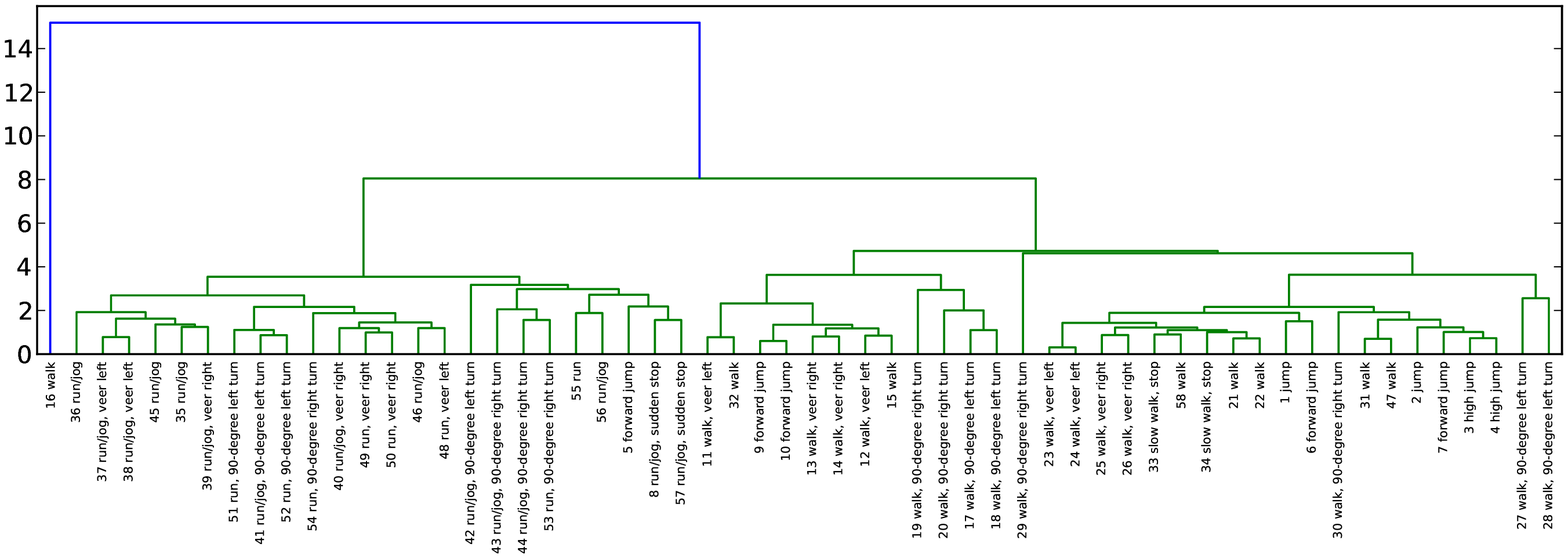}
    \end{minipage}
    \hspace{-1cm}
    %\hfill
    \begin{minipage}[r]{.6\textwidth}
        %\centering
        %\includegraphics[scale=0.8, angle=-90]{cluster16_geodesic.pdf}
        \includegraphics[scale=0.8, angle=-90, trim={0cm 0 1cm 0}, clip]{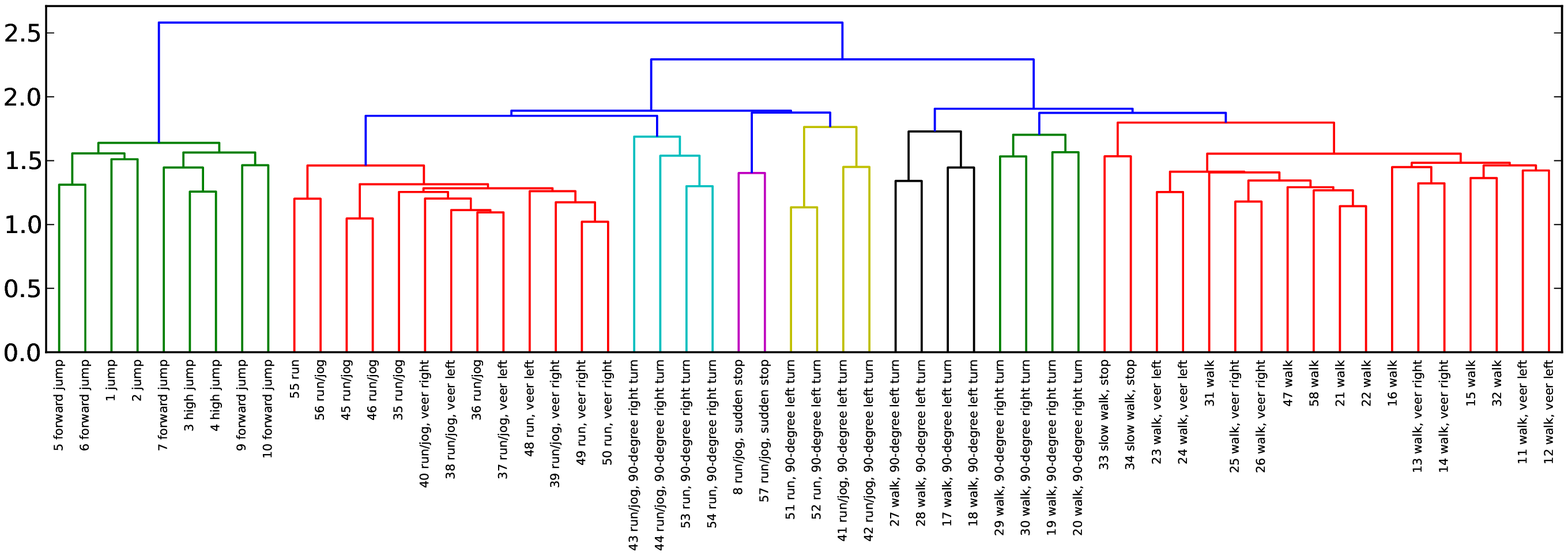}
    \end{minipage}

    \caption{\label{Fig:LinearCluster} Clustering of a variety of walking, running and jumping motions from the CMU database, using normed linear distances on the left side and geodesic distances on the right side. Note particularly the perceived similarities between jumping and walking/running motions in the linear case on the left side, which are not present in the geodesic results on the right side}
\end{figure*}

A hierarchical clustering method was then used on the resulting distance matrix to identify subsets of similar motions.
Some typical results of this can be seen in Figure \ref{Fig:LinearCluster}.
On the left side of that figure, we see a results for the linear case, i.e., where the distance matrix was calculated using linear distances using the $L^2$ norm, whereas geodesic distances were used on the right side.
The same clustering method was used in both cases.
In these tests, we used the space of cyclic animations, $\mathcal{C}^c$, to perform the underlying geodesic computations.
All animations were projected into that space using the projection operator $\mathcal{P}^c$ (\ref{Eq:ProjectionOp}).
It is interesting to note, that this works surprisingly well for non-cyclic animations as well, as can be seen in Figure \ref{Fig:LinearCluster}.

We can immediately see that the geodesic distance measure results in more discrete and clearly delineated clusters of similar motions, whereas the linear distance results in some non-optimal matches, such as a perceived similarity of jumping and walking motions.
Additional tests with other clustering algorithms provided similar results to those in Figure \ref{Fig:LinearCluster}.
Tests using the geodesic distance on $\mathcal{C}^c$, i.e., without accounting for reparametrization, produced slightly worse results than using geodesic distances on motion space $\mathcal{M}$, but still avoided erroneous outliers such as the jumping motions and the animation ``16 walk'' in the normed linear distances test.

%Conclusion.tex
%Author: Markus Eslitzbichler
%Paper: Modelling Character Motions on Infinite-Dimensional Manifolds
%Submitted to: The Visual Computer

\section{Conclusion}\label{Sec:Conclusion}
In this paper, we have proposed the use of a mathematical framework that treats animations as points on infinite dimensional Riemannian manifolds.
This formalism allows us to calculate geodesic paths between animations on those manifolds.
The lengths of such paths can be used to quantify the difference between animations from an energy minimizing point of view.

As we have demonstrated in our examples, this has very widespread applications, while maintaining an elegant simplicity and generality.

In Section \ref{Sec:CycleClosing}, we used only the formal definition of manifolds of open and closed curves\footnote{Corresponding to general and cyclic animations respectively.} to derive a simple and efficient method that can be used to improve the periodicity of animations. In order to blend between two different animations in Sect.~\ref{Sec:Blending}, we then used the concept of geodesic paths on the previously defined manifolds to generate an entire range of blended animations.
Both of these applications use only very general mathematical concepts, but are nevertheless able to produce visually pleasing results.

Finally, in Sect.~\ref{Sec:Classification}, we demonstrated how the manifold modeling approach can be used effectively to analyze and cluster animations.
This has many possible applications, not only in the field of computer graphics, but also for example in biomedical applications, where it could be used in topics such as biometrics and gait recognition.

We believe that the differential geometric approach to processing computer animation, with a focus on energy minimization, has many interesting features to offer and are confident that many applications can profit from the manifold formalism.

Future work on this topic could include handling constraints within the manifold formalism - both joint-angle constraints, as well as physical world space constraints.
The right choice of basis, compatible with the manifold structure, should also be investigated further.
Also the annotation of animation curves with information, e.g., ``left foot is down'', to facilitate the blending of more widely different animations deserves further study.
%Acknowledgements.tex
%Author: Markus Eslitzbichler
%Paper: Modelling Character Motions on Infinite-Dimensional Manifolds
%Submitted to: The Visual Computer

\section{Acknowledgments}\label{Sec:Acknowledgements}
The author would like to thank Elena Celledoni and Markus Grasmair for valuable discussions and feedback.

This research was supported in part by the GeNuIn Applications project grant from the Research Council of Norway.

The data used in this project was obtained from \emph{mocap.cs.cmu.edu}.
The database was created with funding from NSF EIA-0196217.

\bibliographystyle{spmpsci}      % mathematics and physical sciences
\bibliography{AnimationPaper}

\begin{thebibliography}{10}
\providecommand{\url}[1]{{#1}}
\providecommand{\urlprefix}{URL }
\expandafter\ifx\csname urlstyle\endcsname\relax
  \providecommand{\doi}[1]{DOI~\discretionary{}{}{}#1}\else
  \providecommand{\doi}{DOI~\discretionary{}{}{}\begingroup
  \urlstyle{rm}\Url}\fi

\bibitem{abdelkader_silhouette-based_2011}
Abdelkader, M.F., Abd-Almageed, W., Srivastava, A., Chellappa, R.:
  Silhouette-based gesture and action recognition via modeling trajectories on
  riemannian shape manifolds.
\newblock Computer Vision and Image Understanding \textbf{115}(3), 439 -- 455
  (2011).
\newblock \doi{http://dx.doi.org/10.1016/j.cviu.2010.10.006}.
\newblock
  \urlprefix\url{http://www.sciencedirect.com/science/article/pii/S1077314210002377}

\bibitem{bauer_overview_2014}
Bauer, M., Bruveris, M., Michor, P.: Overview of the geometries of shape spaces
  and diffeomorphism groups.
\newblock Journal of Mathematical Imaging and Vision pp. 1--38 (2014).
\newblock \doi{10.1007/s10851-013-0490-z}.
\newblock \urlprefix\url{http://dx.doi.org/10.1007/s10851-013-0490-z}

\bibitem{bauer_sobolev_2011}
Bauer, M., Harms, P., Michor, P.W.: Sobolev metrics on shape space of surfaces.
\newblock Journal of Geometric Mechanics \textbf{3}(4), 389 -- 438 (2011)

\bibitem{bruderlin_motion_1995}
Bruderlin, A., Williams, L.: Motion signal processing.
\newblock In: Proceedings of the 22nd annual conference on Computer graphics
  and interactive techniques, p. 97–104. {ACM} (1995)

\bibitem{carnegie-melllon_carnegie-mellon_2003}
{Carnegie-Mellon}: Carnegie-mellon mocap database. (2003).
\newblock \urlprefix\url{http://mocap.cs.cmu.edu/}

\bibitem{gonzalez_castro_cyclic_2010}
Gonz\'alez~Castro, G., Athanasopoulos, M., Ugail, H.: Cyclic animation using
  partial differential equations.
\newblock The Visual Computer \textbf{26}(5), 325--338 (2010).
\newblock \doi{10.1007/s00371-010-0422-5}.
\newblock \urlprefix\url{http://dx.doi.org/10.1007/s00371-010-0422-5}

\bibitem{kilian_geometric_2007}
Kilian, M., Mitra, N.J., Pottmann, H.: Geometric modeling in shape space.
\newblock {ACM} Transactions on Graphics ({SIGGRAPH)} \textbf{26}(3), \#64,
  1--8 (2007)

\bibitem{klassen_analysis_2004}
Klassen, E., Srivastava, A., Mio, M., Joshi, S.: Analysis of planar shapes
  using geodesic paths on shape spaces.
\newblock {IEEE} Transactions on Pattern Analysis and Machine Intelligence
  \textbf{26}(3), 372 --383 (2004).
\newblock \doi{10.1109/TPAMI.2004.1262333}

\bibitem{kovar_flexible_2003}
Kovar, L., Gleicher, M.: Flexible automatic motion blending with registration
  curves.
\newblock In: Proceedings of the 2003 {ACM} {SIGGRAPH/Eurographics} Symposium
  on Computer Animation, {SCA} '03, p. 214–224. Eurographics Association,
  Aire-la-Ville, Switzerland, Switzerland (2003).
\newblock \urlprefix\url{http://dl.acm.org/citation.cfm?id=846276.846307}

\bibitem{kovar_motion_2002}
Kovar, L., Gleicher, M., Pighin, F.: Motion graphs.
\newblock {ACM} Trans. Graph. \textbf{21}(3), 473–482 (2002).
\newblock \doi{10.1145/566654.566605}.
\newblock \urlprefix\url{http://doi.acm.org/10.1145/566654.566605}

\bibitem{kurtek_novel_2010}
Kurtek, S., Klassen, E., Ding, Z., Srivastava, A.: A novel riemannian framework
  for shape analysis of {3D} objects.
\newblock In: 2010 {IEEE} Conference on Computer Vision and Pattern Recognition
  ({CVPR)}, pp. 1625 --1632 (2010).
\newblock \doi{10.1109/CVPR.2010.5539778}

\bibitem{kurtek_statistical_2012}
Kurtek, S., Srivastava, A., Klassen, E., Ding, Z.: Statistical modeling of
  curves using shapes and related features.
\newblock Journal of the American Statistical Association \textbf{107}(499),
  1152–1165 (2012)

\bibitem{liu_mathematical_2011}
Liu, W., Srivastava, A., Zhang, J.: A mathematical framework for protein
  structure comparison.
\newblock {PLoS} Comput Biol \textbf{7}(2), e1001,075 (2011).
\newblock \doi{10.1371/journal.pcbi.1001075}.
\newblock \urlprefix\url{http://dx.doi.org/10.1371%2Fjournal.pcbi.1001075}

\bibitem{milnor_morse_1963}
Milnor, J.W.: Morse Theory.({AM-51)}, vol.~51.
\newblock Princeton university press (1963)

\bibitem{ormoneit_representing_2005}
Ormoneit, D., Black, M.J., Hastie, T., Kjellstr\"om, H.: Representing cyclic
  human motion using functional analysis.
\newblock Image Vision Comput. \textbf{23}(14), 1264–1276 (2005).
\newblock \doi{10.1016/j.imavis.2005.09.004}.
\newblock \urlprefix\url{http://dx.doi.org/10.1016/j.imavis.2005.09.004}

\bibitem{palais_morse_1963}
Palais, R.S.: Morse theory on hilbert manifolds.
\newblock Topology \textbf{2}(4), 299--340 (1963).
\newblock \doi{10.1016/0040-9383(63)90013-2}.
\newblock
  \urlprefix\url{http://www.sciencedirect.com/science/article/pii/0040938363900132}

\bibitem{pejsa_state_2010}
Pejsa, T., Pandzic, I.: State of the art in example-based motion synthesis for
  virtual characters in interactive applications.
\newblock Computer Graphics Forum \textbf{29}(1), 202–226 (2010).
\newblock \doi{10.1111/j.1467-8659.2009.01591.x}.
\newblock \urlprefix\url{http://dx.doi.org/10.1111/j.1467-8659.2009.01591.x}

\bibitem{shoemake_animating_1985}
Shoemake, K.: Animating rotation with quaternion curves.
\newblock {SIGGRAPH} Comput. Graph. \textbf{19}(3), 245–254 (1985).
\newblock \doi{10.1145/325165.325242}.
\newblock \urlprefix\url{http://doi.acm.org/10.1145/325165.325242}

\bibitem{srivastava_shape_2011}
Srivastava, A., Klassen, E., Joshi, S., Jermyn, I.: Shape analysis of elastic
  curves in euclidean spaces.
\newblock {IEEE} Transactions on Pattern Analysis and Machine Intelligence
  \textbf{33}(7), 1415 --1428 (2011).
\newblock \doi{10.1109/TPAMI.2010.184}

\bibitem{srivastava_advances_2012}
Srivastava, A., Turaga, P., Kurtek, S.: On advances in differential-geometric
  approaches for {2D} and {3D} shape analyses and activity recognition.
\newblock Image and Vision Computing \textbf{30}(6–7), 398--416 (2012).
\newblock \doi{10.1016/j.imavis.2012.03.006}.
\newblock
  \urlprefix\url{http://www.sciencedirect.com/science/article/pii/S0262885612000492}

\bibitem{wei_liu_riemannian_????}
{Wei Liu}: A riemannian framework for annotated curves analysis.
\newblock Ph.D. thesis, The Florida State University.
\newblock \urlprefix\url{http://diginole.lib.fsu.edu/etd/4997}

\bibitem{younes_spaces_2012}
Younes, L.: Spaces and manifolds of shapes in computer vision: An overview.
\newblock Image and Vision Computing \textbf{30}(6–7), 389--397 (2012).
\newblock \doi{10.1016/j.imavis.2011.09.009}.
\newblock
  \urlprefix\url{http://www.sciencedirect.com/science/article/pii/S0262885611001028}

\end{thebibliography}
\end{document}